\documentclass[12pt,twoside]{amsart}
\usepackage{amsmath}
\usepackage{amsthm}
\usepackage{amsfonts}
\usepackage{amssymb}
\usepackage{latexsym}
\usepackage[all]{xy}

\date{}
\allowdisplaybreaks[4] \footskip=15pt
\renewcommand{\uppercasenonmath}[1]{}

\numberwithin{equation}{section} \theoremstyle{plain}
\newtheorem*{thm*}{Main Theorem}
\newtheorem{thm}{Theorem}[section]
\newtheorem{cor}[thm]{Corollary}
\newtheorem*{cor*}{Corollary}
\newtheorem{lem}[thm]{Lemma}
\newtheorem*{lem*}{Lemma}

\newtheorem*{prop*}{Proposition}

\newtheorem*{rem*}{Remark}
\newtheorem{exa}[thm]{Example}
\newtheorem*{exa*}{Example}
\newtheorem{df}[thm]{Definition}
\newtheorem*{df*}{Definition}

\newtheorem*{conj*}{Conjecture}
\newtheorem*{ack*}{ACKNOWLEDGEMENTS}

\newcommand{\pf}{\noindent\begin {proof}}
\newcommand{\epf}{\end{proof}}

\newcommand{\Ext}{\mbox{\rm Ext}}

\newcommand{\Hom}{\mbox{\rm Hom}}

\newcommand{\im}{\mbox{\rm im}}
\newcommand{\coker}{\mbox{\rm coker}}
\begin{document}
\begin{center}
{\bf \large Gorenstein projective, injective and flat modules over trivial ring extensions}

\vspace{0.5cm} Lixin Mao\\
School of Mathematics and Physics, Nanjing Institute of
Technology,\\ Nanjing 211167, China\\
E-mail: maolx2@hotmail.com \\
\end{center}

\bigskip
\centerline { \bf  Abstract}
 \bigskip
\leftskip10truemm \rightskip10truemm
 \noindent
We introduce the concepts of generalized compatible and cocompatible bimodules in order to characterize  Gorenstein projective, injective and flat modules over trivial ring extensions. Let $R\ltimes M$ be a trivial  extension of a ring $R$ by an $R$-$R$-bimodule $M$ such that $M$ is a generalized compatible $R$-$R$-bimodule and $\textbf{Z}(R)$ is a generalized  compatible $R\ltimes M$-$R\ltimes M$-bimodule. We prove that $(X,\alpha)$ is a Gorenstein projective left $R\ltimes M$-module if and only if the sequence $M\otimes_R M\otimes_R X\stackrel{M\otimes\alpha}\rightarrow M\otimes_R X\stackrel{\alpha}\rightarrow X$ is exact and $\coker(\alpha)$ is a Gorenstein projective left $R$-module. Analogously, we explicitly characterize Gorenstein injective and flat modules over trivial ring extensions. As an application, we describe Gorenstein projective, injective and flat modules over  Morita context rings with  zero bimodule homomorphisms.\\
\vbox to 0.3cm{}

\leftskip0truemm \rightskip0truemm
\bigskip
\section { \bf Introduction}
\bigskip
Let $R$ be an associative ring and $M$ be an $R$-$R$-bimodule. Then the Cartesian product $R \times M$, with the natural addition and the multiplication given by $(r_{1}, m_{1})(r_{2}, m_{2}) = (r_{1}r_{2}, r_{1}m_{2}+m_{1}r_{2})$, becomes a ring. This ring is
called the \emph{trivial  extension} of the ring $R$ by the bimodule $M$, and denoted
by $R\ltimes M$. When $R$ is a commutative ring, Nagata also called this construction an \emph{idealization} in \cite{N}. The notion of trivial extension of a ring by a bimodule is an important extension of rings and has played a crucial role in ring theory and homological algebra \cite{DMZ, FGR, IM,MY, N, PR, R1}. Fossum, Griffith and Reiten studied the categorical aspect and homological properties of trivial ring extensions \cite{FGR}. Palm\'{e}r and Roos gave some explicit formulae for the global homological dimensions of trivial ring extensions \cite{PR}. Recently, Minamoto and Yamaura furthermore investigated homological dimension formulas for trivial extension algebras \cite{MY}. Dumitrescu, Mahdou and  Zahir studied the radical factorization for trivial extensions \cite{DMZ}.

The origin of Gorenstein homological algebra may date back to 1960s
when Auslander and Bridger introduced the concept of G-dimension for
finitely generated modules over a two-sided Noetherian ring
\cite{AB}. In 1990s, Enochs, Jenda and Torrecillas extended the
ideas of Auslander and Bridger and introduced the concepts of
Gorenstein projective, injective and flat modules over arbitrary
rings, and then developed Gorenstein homological algebra \cite{EJ1, EJ, EJT}. When $R$ is a  commutative Noetherian ring and $M$ is a semi-dualizing $R$-module,  Holm and J{\o}rgensen investigated Gorenstein projective, injective and flat modules over a trivial ring extension $R\ltimes M$ \cite{HJ} and they also obtained more general results in categories of quiver representations \cite{HJ1}. Mahdou and Ouarghi studied the Gorenstein dimensions of a trivial ring extension $R\ltimes M$ when $R$ is a  commutative ring \cite{MO}.

In the present paper, we further consider the transfer of Gorenstein properties between an arbitrary ring $R$ (not necessarily commutative) and its trivial ring extension $R\ltimes M$. The contents of this paper are arranged as follows.

In Section 2, we recall some basic concepts and results on trivial extensions.

In Section 3, we first introduce the concept of generalized compatible  bimodules to describe Gorenstein projective $R\ltimes M$-modules. Let $M$ be a  generalized compatible $R$-$R$-bimodule and $\textbf{Z}(R)=(R,0)$ be a generalized  compatible $R\ltimes M$-$R\ltimes M$-bimodule. Then we prove that $(X,\alpha)$ is a Gorenstein projective left $R\ltimes M$-module if and only if the sequence $M\otimes_R M\otimes_R X\stackrel{M\otimes\alpha}\rightarrow M\otimes_R X\stackrel{\alpha}\rightarrow X$ is exact and $\coker(\alpha)$ is a Gorenstein projective left $R$-module.

Section 4 is devoted to Gorenstein injective and flat $R\ltimes M$-modules. We first introduce the concept of generalized cocompatible  bimodules. Let $M$ be a  generalized cocompatible $R$-$R$-bimodule and $\textbf{Z}(R)=(R,0)$ be a generalized  cocompatible $R\ltimes M$-$R\ltimes M$-bimodule. Then we prove that: (1) $[Y,\beta]$  is a Gorenstein injective left $R\ltimes M$-module if and only if the sequence  $Y\stackrel{\beta}\rightarrow \Hom_R(M, Y)\stackrel{\beta_{*}}\rightarrow \Hom_R(M,\Hom_R(M, Y))$ is exact and $\ker(\beta)$ is a Gorenstein injective left $R$-module; (2) if $R\ltimes M$ is a left coherent ring, then $(X,\alpha)$ is a Gorenstein flat right $R\ltimes M$-module if and only if the sequence $X\otimes_R M\otimes_R M \stackrel{\alpha\otimes M}\rightarrow X\otimes_R M \stackrel{\alpha}\rightarrow X$  is exact and $\coker(\alpha)$ is a Gorenstein flat right $R$-module.

In Section 5, we apply the foregoing results to one important special case of trivial ring extensions, i.e., Morita context rings with  zero bimodule homomorphisms. In \cite{GP}, Gao and Psaroudakis investigated how to construct Gorenstein projective modules over a Morita context ring which is an Artin algebra and has  zero bimodule homomorphisms. We extend the result to a general setting and describe Gorenstein projective, injective and flat modules over Morita context rings with  zero bimodule homomorphisms.
\bigskip
\section {\bf Preliminaries and notations}
\bigskip
Throughout this paper, all rings are nonzero associative rings with identity and all modules are unitary. For a ring $R$, $_RX$ (resp. $X_{R}$) denotes a left (resp. right) $R$-module. We write $R$-Mod (resp. Mod-$R)$  for the category of left (resp. right) $R$-modules. For an $R$-module $X$, the character module $\Hom_{\mathbb{Z}}(X,\mathbb{Q}/\mathbb{Z})$ of $X$  is denoted by $X^{+}$, Add($X$) (resp. Prod($X$))  denotes the class of $R$-modules isomorphic to direct summands of direct sums  (resp. direct products) of copies of $X$, $pd(X)$, $id(X)$ and $fd(X)$ denote the projective, injective and flat dimensions of $X$, respectively.

Recall from \cite{FGR} that the classical \emph{right trivial extension} of an abelian category $\textbf{\underline{\underline{A}}}$ by an additive  endofunctor $\textbf{F}$, denoted by $\textbf{\underline{\underline{A}}}\ltimes \textbf{F}$, is a category whose objects are couples $(X, f)$ with $X\in$ Ob$(\textbf{\underline{\underline{A}}})$ and $f: \textbf{F}(X) \rightarrow X$ such that $f\cdot \textbf{F}(f) = 0$ and a morphism $\gamma: (X,\alpha) \rightarrow (Y,\beta)$ is a morphism $\gamma: X\rightarrow Y$ in $\textbf{\underline{\underline{A}}}$ such that the following diagram $$\xymatrix{\textbf{F}(X)\ar[d]_{\alpha}\ar[r]^{\textbf{F}(\gamma)}&\textbf{F}(Y)\ar[d]_{\beta}\\
X\ar[r]^{\gamma}&Y}$$ commutes. If $\textbf{F}$ is right exact, then  $\textbf{\underline{\underline{A}}}\ltimes \textbf{F}$ is an abelian category. In this case, a sequence in $\textbf{\underline{\underline{A}}}\ltimes \textbf{F}$ is exact if and only if the sequence of codomains in $\textbf{\underline{\underline{A}}}$ is exact.

Dually, the \emph{left trivial extension} of an abelian category $\textbf{\underline{\underline{A}}}$ by an additive  endofunctor $\textbf{G}$, denoted by $\textbf{G}\rtimes \textbf{\underline{\underline{A}}}$, is a category whose objects are couples $[X, g]$ with $X \in$ Ob$(\textbf{\underline{\underline{A}}})$ and $g: X \rightarrow \textbf{G}(X)$ such that $\textbf{G}(g)\cdot g = 0$ and a morphism $\gamma: [X,\alpha] \rightarrow [Y,\beta]$ is a morphism $\gamma:X\rightarrow Y$ in $\textbf{\underline{\underline{A}}}$ such that the following diagram $$\xymatrix{\textbf{G}(X)\ar[r]^{\textbf{G}(\gamma)}&\textbf{G}(Y)\\X\ar[u]_{\alpha}\ar[r]^{\gamma}&Y\ar[u]_{\beta}}$$ commutes. If $\textbf{G}$ is left exact, then  $\textbf{G}\rtimes \textbf{\underline{\underline{A}}}$ is an abelian category. In this case, a sequence in $\textbf{G}\rtimes \textbf{\underline{\underline{A}}}$ is exact if and only if the sequence of domains in $\textbf{\underline{\underline{A}}}$ is exact.

For a right exact  endofunctor $\textbf{F}: \textbf{\underline{\underline{A}}}\rightarrow \textbf{\underline{\underline{A}}}$ and a left exact  endofunctor $\textbf{G}: \textbf{\underline{\underline{A}}}\rightarrow \textbf{\underline{\underline{A}}}$, there are some important functors as follows.

The functor $\textbf{T}: \textbf{\underline{\underline{A}}}\rightarrow \textbf{\underline{\underline{A}}}\ltimes \textbf{F}$ is given, for every object $X\in \textbf{\underline{\underline{A}}}$, by $\textbf{T}(X) =(X\oplus \textbf{F}(X), \mu)$ with $\mu=\biggl(\begin{matrix} 0&0\\1&0 \end{matrix}\biggr): \textbf{F}(X)\oplus \textbf{F}^{2}(X)\rightarrow X\oplus \textbf{F}(X)$ and for morphisms by $\textbf{T}(\alpha)=\biggl(\begin{matrix} \alpha&0\\0&\textbf{F}(\alpha)\end{matrix}\biggr)$.

The functor $\textbf{U}:\textbf{\underline{\underline{A}}}\ltimes \textbf{F}\rightarrow \textbf{\underline{\underline{A}}}$ is given, for every object $(X,f)\in \textbf{\underline{\underline{A}}}\ltimes \textbf{F}$, by $\textbf{U}(X,f) =X$ and for morphisms by $\textbf{U}(\alpha)=\alpha$.

The functor $\textbf{Z}: \textbf{\underline{\underline{A}}}\rightarrow \textbf{\underline{\underline{A}}}\ltimes \textbf{F}$ is given, for every object  $X\in \textbf{\underline{\underline{A}}}$, by $\textbf{Z}(X)=(X,0)$ and for morphisms by $\textbf{Z}(\alpha)=\alpha$.

The functor $\textbf{C}:\textbf{\underline{\underline{A}}}\ltimes \textbf{F}\rightarrow \textbf{\underline{\underline{A}}}$ is given, for every object $(X,f)\in \textbf{\underline{\underline{A}}}\ltimes \textbf{F}$, by $\textbf{C}(X,f) =\coker(f)$ and for morphisms by $\textbf{C}(\alpha)=$ the induced morphism.

The functor $\textbf{H}: \textbf{\underline{\underline{A}}}\rightarrow \textbf{G}\rtimes\textbf{\underline{\underline{A}}}$ is given, for every object  $X\in \textbf{\underline{\underline{A}}}$, by $\textbf{H}(X)=[\textbf{G}(X)\oplus X,\vartheta]$ with $\vartheta=\biggl(\begin{matrix} 0&0\\1&0 \end{matrix}\biggr): \textbf{G}(X)\oplus X\rightarrow \textbf{G}^{2}(X)\oplus \textbf{G}(X)$ and for morphisms by $\textbf{H}(\beta)=\biggl(\begin{matrix} \textbf{G}(\beta)&0\\0&\beta\end{matrix}\biggr)$.

The functor $\textbf{U}:\textbf{G}\rtimes\textbf{\underline{\underline{A}}}\rightarrow \textbf{\underline{\underline{A}}}$ is given, for every object $[X,g]\in \textbf{G}\rtimes\textbf{\underline{\underline{A}}}$, by $\textbf{U}[X,g] =X$ and for morphisms by $\textbf{U}(\alpha)=\alpha$.

The functor $\textbf{Z}: \textbf{\underline{\underline{A}}}\rightarrow \textbf{G}\rtimes\textbf{\underline{\underline{A}}}$ is given, for every object  $X\in \textbf{\underline{\underline{A}}}$, by $\textbf{Z}(X)=[X,0]$ and for morphisms by $\textbf{Z}(\alpha)=\alpha$.

The functor $\textbf{K}: \textbf{G}\rtimes\textbf{\underline{\underline{A}}}\rightarrow \textbf{\underline{\underline{A}}}$ is given, for every object $[X,g]\in \textbf{G}\rtimes\textbf{\underline{\underline{A}}}$, by $\textbf{K}[X,g] =\ker(g)$ and for morphisms by $\textbf{K}(\alpha)=$ the induced morphism.

We note that the functors $\textbf{T}$ and $\textbf{C}$ are right exact, $\textbf{H}$ and $\textbf{K}$ are left exact, $\textbf{U}$ and $\textbf{Z}$ are exact. There exist important pairs of adjoint functors $(\textbf{T}, \textbf{U})$, $(\textbf{C}, \textbf{Z})$, $(\textbf{Z}, \textbf{K})$ and $(\textbf{U}, \textbf{H})$ for which $\textbf{C}\textbf{T}= id_{\textbf{\underline{\underline{A}}}}$, $\textbf{U}\textbf{Z}= id_{\textbf{\underline{\underline{A}}}}$, $\textbf{K}\textbf{H}= id_{\textbf{\underline{\underline{A}}}}$. $$\xymatrix{&\textbf{\underline{\underline{A}}}\ar@<6pt>[r]^{\textbf{T}}&\textbf{\underline{\underline{A}}}\ltimes \textbf{F}\ar@<6pt>[l]^{\textbf{U}}\ar@<6pt>[r]^{\textbf{C}}&\textbf{\underline{\underline{A}}}\ar@<6pt>[l]^{\textbf{Z}}},\xymatrix{&\textbf{\underline{\underline{A}}}\ar@<6pt>[r]^{\textbf{Z}}&\textbf{G}\rtimes\textbf{\underline{\underline{A}}}\ar@<6pt>[l]^{\textbf{K}}\ar@<6pt>[r]^{\textbf{U}}
&\textbf{\underline{\underline{A}}}\ar@<6pt>[l]^{\textbf{H}}}.$$

It is known that, when $\textbf{\underline{\underline{A}}}$ is the category of left $R$-modules, $M$ is an $R$-$R$-bimodule, $\textbf{F}=M\otimes_{R}-$ and $\textbf{G}=\Hom_{R}(M,-)$, both $\textbf{\underline{\underline{A}}}\ltimes \textbf{F}$ and $\textbf{G}\rtimes \textbf{\underline{\underline{A}}}$ are isomorphic to the category of left modules over $R\ltimes M$. We will identify $R\ltimes M$-Mod with $\textbf{\underline{\underline{A}}}\ltimes \textbf{F}$ and $\textbf{G}\rtimes \textbf{\underline{\underline{A}}}$ in what follows.
\begin{lem} \label{lem: 2.1}\cite[Corollary 1.6 (c) and (d)]{FGR} Let $R\ltimes M$ be a trivial  extension of a ring $R$ by an $R$-$R$-bimodule $M$. Then
\begin{enumerate}\item $(X,\alpha)$ is a projective left $R\ltimes M$-module if and only if $(X,\alpha)\cong\textbf{T}(P)$ for a projective left $R$-module $P$.\item $[Y,\beta]$ is an injective left $R\ltimes M$-module if and only if $[Y,\beta]\cong\textbf{H}(E)$ for an injective left $R$-module $E$.
\end{enumerate}
\end{lem}
\begin{lem} \label{lem: 2.2} Let $(X,\alpha)$ and $[Y,\beta]$ be left $R\ltimes M$-modules, $\rho: X\rightarrow \coker(\alpha)$ be the canonical epimorphism, $\iota: \ker(\beta)\rightarrow Y$ be the inclusion.
\begin{enumerate}\item $\textbf{Z}(W)\otimes_{R\ltimes M}(X,\alpha)\cong W\otimes_{R}\coker(\alpha)$ for any right $R$-module $W$.\item
$\Hom_{R\ltimes M}(\textbf{Z}(X),[Y,\beta])\cong \Hom_{R}(X,\ker(\beta))$.\item There is an exact sequence $0\rightarrow \textbf{Z}(\im(\alpha))\rightarrow (X,\alpha)\rightarrow \textbf{Z}(\coker(\alpha))\rightarrow 0$.\item There is an exact sequence $0\rightarrow \textbf{Z}(\ker(\beta))\rightarrow [Y,\beta]\rightarrow \textbf{Z}(\im(\beta))\rightarrow 0$.
\item There is $\delta: M\otimes_{R}\coker(\alpha)\rightarrow X$ such that $\delta(M\otimes\rho)=\alpha$.\item There is $\gamma: Y\rightarrow \Hom_{R}(M,\ker(\beta))$ such that $\iota_{*}\gamma=\beta$, where the morphism $\iota_{_{*}}: \Hom_{R}(M,\ker(\beta))\rightarrow \Hom_R(M, Y)$ is induced by $\iota$.
\end{enumerate}
\end{lem}
\begin{proof}(1) Note that $\textbf{Z}(W)\cong$ $ W\otimes_{R}R_{R\ltimes M}$ and $_{R}R\otimes_{R\ltimes M}(X,\alpha)\cong \coker(\alpha)$. So
$$\textbf{Z}(W)\otimes_{R\ltimes M}(X,\alpha)\cong(W\otimes_{R}R_{R\ltimes M}) \otimes_{R\ltimes M} (X,\alpha)\cong W\otimes_{R}(R\otimes_{R\ltimes M}(X,\alpha))\cong W\otimes_{R}\coker(\alpha).$$

(2) is clear by the adjoint pair $(\textbf{Z}, \textbf{K})$.

(3) follows from the following commutative diagram: $$\xymatrix{& M\otimes_R \im(\alpha)\ar[r]\ar[d]_{0}&M\otimes_R X\ar[d]_{\alpha}\ar[r]&M\otimes_R\coker(\alpha)\ar[d]_{0}\\
0\ar[r]&\im(\alpha)\ar[r]&X\ar[r]&\coker(\alpha)\ar[r]&0.}$$

(4) follows from the following commutative diagram: $$\xymatrix{0\ar[r]&\ker(\beta)\ar[r]\ar[d]_{0}&Y\ar[d]_{\beta}\ar[r]&\im(\beta)\ar[d]_{0}\ar[r]&0\\
&\Hom_ {R}(M,\ker(\beta))\ar[r]&\Hom_{R}(M,Y)\ar[r]&\Hom_{R}(M,\im(\beta)).}$$

(5) The exact seqence $M\otimes_R X\stackrel{\alpha}\rightarrow X\stackrel{\rho}\rightarrow \coker(\alpha)\rightarrow 0$ induces the exact sequence $M\otimes_R M\otimes_R X\stackrel{M\otimes\alpha}\rightarrow M\otimes_R X\stackrel{M\otimes\rho}\rightarrow M\otimes_{R}\coker(\alpha)\rightarrow 0$. Since $\alpha (M\otimes\alpha)=0$, there is $\delta: M\otimes_{R}\coker(\alpha)\rightarrow X$ such that $\delta(M\otimes\rho)=\alpha$.

(6) The exact sequence $0\rightarrow \ker(\beta)\stackrel{\iota}\rightarrow Y\stackrel{\beta}\rightarrow \Hom_R(M, Y)$ induces the exact sequence  $0\rightarrow \Hom_{R}(M,\ker(\beta))\stackrel{\iota_{_{*}}}\rightarrow \Hom_R(M, Y)\stackrel{\beta_{*}}\rightarrow \Hom_R(M,\Hom_R(M, Y))$. Since $\beta_{*}\beta=0$, there is $\gamma: Y\rightarrow \Hom_{R}(M,\ker(\beta))$ such that $\iota_{*}\gamma=\beta$.
\end{proof}
\bigskip
\section {\bf Gorenstein projective modules over trivial ring extensions}
\bigskip
Recall that an exact sequence of projective left $R$-modules $$\Xi: \cdots\rightarrow
P^{-1}\stackrel{f^{-1}}\rightarrow P^{0}\stackrel{f^{0}}\rightarrow P^{1}\stackrel{f^{1}}\rightarrow P^{2}\rightarrow \cdots$$ is a \emph{complete projective resolution} in $R$-Mod \cite{H} if $\Hom_{R}(\Xi, Q)$ is also exact for any projective left $R$-module $Q$.

A left $R$-module $X$ is called \emph{Gorenstein projective} \cite{EJ1} if there is a complete projective resolution $\Xi: \cdots\rightarrow
P^{-1}\stackrel{f^{-1}}\rightarrow P^{0}\stackrel{f^{0}}\rightarrow P^{1}\stackrel{f^{1}}\rightarrow P^{2}\rightarrow \cdots$ in $R$-Mod  such that $X\cong\ker(f^{0})$.

In \cite{Z}, the concept of  compatible bimodules was introduced by Zhang in order to describe Gorenstein projective modules over triangular matrix algebras. Note that a triangular matrix algebra is a special trivial ring extension. Here we give the concept of generalized compatible bimodules in order to describe Gorenstein projective modules over  trivial ring extensions.
\begin{df} \label{df: 3.1} {\rm $_{A}N_{B}$ is called a \emph{generalized compatible bimodule} if the following two conditions hold:\\
(1) If $\Xi$ is a complete projective resolution in $B$-Mod,  then $N\otimes_{B}\Xi$ is exact.\\
(2) If $\Delta$ is a complete projective resolution in $A$-Mod, then $\Hom_{A}(\Delta$, Add$(N))$ is exact.}
\end{df}
\begin{exa} \label{exa: 3.2}{\rm (1) If $fd(N_{B})<\infty$ and $pd(_{A}N)<\infty$, then $_{A}N_{B}$ is a generalized compatible bimodule by \cite[Lemma 2.3]{EIT} and \cite[Proposition 2.3]{H}.

(2) If $A$ is a left Noetherian ring, $fd(N_{B})<\infty$ and $id(_{A}N)<\infty$, then $_{A}N_{B}$ is a generalized compatible bimodule by \cite[Lemmas 2.3 and 2.4]{EIT}. }
\end{exa}
\begin{thm} \label{thm: 3.3}Let $M$ be a  generalized compatible $R$-$R$-bimodule and $(X,\alpha)$ be a left $R\ltimes M$-module such that the sequence $M\otimes_R M\otimes_R X\stackrel{M\otimes\alpha}\rightarrow M\otimes_R X\stackrel{\alpha}\rightarrow X$ is exact and $\coker(\alpha)$ is a Gorenstein projective left $R$-module. Then  $(X,\alpha)$ is a Gorenstein projective left $R\ltimes M$-module.
\end{thm}
\begin{proof}There exists a complete projective resolution  in $R$-Mod
$$\Xi: \cdots\rightarrow
P^{-1}\stackrel{f^{-1}}\rightarrow
P^{0}\stackrel{f^{0}}\rightarrow
P^{1}\stackrel{f^{1}}\rightarrow
P^{2}\rightarrow \cdots$$ with
$\coker(\alpha)\cong\ker(f^{0})$. Since $M$ is a  generalized compatible $R$-$R$-bimodule, we get the exact sequence of  left $R$-modules
$$M\otimes_{R}\Xi: \cdots\rightarrow
M\otimes_{R} P^{-1}\stackrel{M\otimes f^{-1}}\rightarrow
M\otimes_{R} P^{0}\stackrel{M\otimes f^{0}}\rightarrow
M\otimes_{R} P^{1}\stackrel{M\otimes f^{1}}\rightarrow
M\otimes_{R} P^{2}\rightarrow \cdots$$ with
$M\otimes_{R}\coker(\alpha)\cong\ker({M\otimes f^{0}})$.

The two exact sequences $$M\otimes_R M\otimes_R X\stackrel{M\otimes\alpha}\rightarrow M\otimes_R X\stackrel{\alpha}\rightarrow X\stackrel{\rho}\rightarrow \coker(\alpha)\rightarrow 0$$  and $$M\otimes_R M\otimes_R X\stackrel{M\otimes\alpha}\rightarrow M\otimes_R X\stackrel{M\otimes\rho}\rightarrow M\otimes_R \coker(\alpha)\rightarrow 0$$ yield the exact sequence $$0\rightarrow M\otimes_R \coker(\alpha)\stackrel{\delta}\rightarrow X\stackrel{\rho}\rightarrow \coker(\alpha)\rightarrow 0$$ with $\delta(M\otimes\rho)=\alpha$ by Lemma \ref{lem: 2.2}(5). Since $M\otimes_{R}P^{i}\in$ Add$(M)$ and $M$ is a  generalized compatible $R$-$R$-bimodule, the complex $\Hom_{R}(\Xi, M\otimes_{R}P^{i})$ is exact. So we have $\Ext^{1}_{R}(\ker(f^{i}), M\otimes_{R}P^{i})=0$. Let $\iota: \coker(\alpha)\rightarrow P^{0}$ be the inclusion and $\pi: P^{-1}\rightarrow\coker(\alpha)$ be the  canonical  epimorphism such that $\iota\pi=f^{-1}$. Then there is $\psi: X\rightarrow M\otimes_{R}P^{0}$ such that $\psi\delta=M\otimes\iota$ and there is $\eta: P^{-1}\rightarrow X$ such that $\rho\eta=\pi$.

Define $\lambda: X\rightarrow P^{0}\oplus (M\otimes_{R}P^{0})$  by $$\lambda(x)=(\iota\rho(x),\psi(x))$$ and define $\xi: P^{-1}\oplus (M\otimes_{R}P^{-1})\rightarrow X$ by $$\xi(x,y)=\eta(x)+\delta(M\otimes\pi)(y).$$ It is easy to check that $\lambda$ is a monomorphism and $\xi$ is an epimorphism. By  the generalized Horseshoe Lemma \cite[Lemma 1.6]{Z}, we get the following commutative diagrams:
$$\xymatrix{&0\ar[d]&&0\ar[d] &0\ar[d] &\\
0 \ar[r] &M\otimes_R \coker(\alpha)\ar[d]^{M\otimes\iota}\ar[rr]^{\delta}&&X\ar@{.>}[d]^{\lambda}\ar[r]^{\rho}&\coker(\alpha)\ar[r]\ar[d]^{\iota}&0\\
0\ar[r] &M\otimes_{R}P^{0}\ar[d]^{^{M\otimes f^{0}}}\ar[rr]&&P^{0}\oplus (M\otimes_{R}P^{0})\ar@{.>}[d]^{g^{0}}
\ar[r]&P^{0}\ar[r]\ar[d]^{f^{0}}&0\\
0\ar[r]&M\otimes_{R}P^{1}\ar[d]^{M\otimes f^{1}}\ar[rr]&&P^{1}\oplus (M\otimes_{R}P^{1})\ar@{.>}[d]^{g^{1}}
\ar[r]&P^{1}\ar[r]\ar[d]^{f^{1}}&0\\
0\ar[r]&M\otimes_{R} P^{2}\ar[d]\ar[rr]&&P^{2}\oplus (M\otimes_{R}P^{2})
\ar@{.>}[d]\ar[r]&P^{2}\ar[d]\ar[r]&0\\
&\vdots &&\vdots &\vdots &}$$
and
$$\xymatrix{&\vdots\ar[d]
&&\vdots\ar@{.>}[d] &\vdots\ar[d] &\\
0\ar[r] &M\otimes_{R}P^{-2}\ar[d]^{^{M\otimes f^{-2}}}\ar[rr]&&P^{-2}\oplus (M\otimes_{R}P^{-2})\ar@{.>}[d]^{g^{-2}}
\ar[r]&P^{-2}\ar[r]\ar[d]^{f^{-2}}&0\\0
\ar[r]&M\otimes_{R}P^{-1}\ar[d]^{M\otimes \pi}\ar[rr]&&P^{-1}\oplus (M\otimes_{R}P^{-1})\ar@{.>}[d]^{\xi}
\ar[r]&P^{-1}\ar[r]\ar[d]^{\pi}&0\\0\ar[r] &M\otimes_R \coker(\alpha)\ar[d]\ar[rr]^{\delta}&&X\ar[d]
\ar[r]^{\rho}&\coker(\alpha)\ar[r]\ar[d]&0\\
&0&&0 &0 &}$$ with exact rows and columns and $g^{i}=\biggl(\begin{matrix} f^{i}&0\\\sigma^{i}&M\otimes f^{i}\end{matrix}\biggr)$.

Let $g^{-1}=\lambda\xi$. Then  we get the exact sequence of  left $R$-modules
$$\cdots\rightarrow P^{-1}\oplus (M\otimes_{R} P^{-1})\stackrel{g^{-1}}\rightarrow P^{0}\oplus
(M\otimes_{R} P^{0})\stackrel{g^{0}}\rightarrow P^{1}\oplus
(M\otimes_{R} P^{1})\stackrel{g^{1}}\rightarrow P^{2}\oplus
(M\otimes_{R} P^{2})\rightarrow \cdots$$  with $X\cong\ker(g^{0})$. The following commutative diagram
$${\Small \xymatrix{&M\otimes_R X\ar[d]_{\alpha}\ar[rr]^{M\otimes\lambda}&&M\otimes_R(P^{0}\oplus
(M\otimes_{R} P^{0}))\ar[d]^{\mu_{0}}\ar[r]^{M\otimes g^{0}}&M\otimes_R(P^{1}\oplus
(M\otimes_{R} P^{1}))\ar[d]^{\mu_{1}}\ar[r]&\cdots\\
0\ar[r]&X\ar[rr]^{\lambda}&&P^{0}\oplus
(M\otimes_{R} P^{0})\ar[r]^{g^{0}}&P^{1}\oplus
(M\otimes_{R} P^{1})\ar[r]&\cdots}}$$ implies that the sequence $0\rightarrow(X,\alpha)\stackrel{\lambda}\rightarrow\textbf{T}(P^{0})\stackrel{g^{0}}\rightarrow \textbf{T}(P^{1})\rightarrow\cdots$ is exact.

The following commutative diagram
$${\Small\xymatrix{\cdots\ar[r]&M\otimes_R(P^{-2}\oplus (M\otimes_{R} P^{-2}))\ar[d]^{\mu_{-2}}\ar[r]^{M\otimes g^{-2}}&M\otimes_R(P^{-1}\oplus
(M\otimes_{R} P^{-1}))\ar[d]^{\mu_{-1}}\ar[rr]^{M\otimes\xi}&&M\otimes_R X\ar[d]_{\alpha}\\
\cdots\ar[r]&P^{-2}\oplus (M\otimes_{R} P^{-2})\ar[r]^{g^{-2}}&P^{-1}\oplus (M\otimes_{R} P^{-1})\ar[rr]^{\xi}&&X\ar[r]&0}}$$ implies that the sequence $\cdots\rightarrow\textbf{T}(P^{-2})\stackrel{g^{-2}}\rightarrow \textbf{T}(P^{-1})\stackrel{\xi}\rightarrow(X,\alpha)\rightarrow 0$ is exact.

By Lemma \ref{lem: 2.1}(1),  each $\textbf{T}(P^{i})$ is projective. So there is an
exact sequence of projective left $R\ltimes M$-modules $$\Delta: \cdots\rightarrow \textbf{T}(P^{-1})\stackrel{g^{-1}}\rightarrow
\textbf{T}(P^{0})\stackrel{g^{0}}\rightarrow \textbf{T}(P^{1})\stackrel{g^{1}}\rightarrow \textbf{T}(P^{2})\rightarrow \cdots$$ with $(X,\alpha)\cong\ker(g^{0})$.

By Lemma \ref{lem: 2.1}(1), any projective left $R\ltimes M$-module is isomorphic to $\textbf{T}(Q)$ with $Q$  a projective left $R$-module. By Lemma \ref{lem: 2.2}(3), there is an exact sequence in $R\ltimes M$-Mod  $$0\rightarrow \textbf{Z}(M\otimes_{R} Q)\rightarrow \textbf{T}(Q)\rightarrow \textbf{Z}(Q)\rightarrow 0.$$ So we get the exact sequence of complexes $$0\rightarrow \Hom_{R\ltimes M}(\Delta,\textbf{Z}(M\otimes_{R}Q))\rightarrow \Hom_{R\ltimes M}(\Delta,\textbf{T}(Q))\rightarrow \Hom_{R\ltimes M}(\Delta,\textbf{Z}(Q))\rightarrow 0.$$
Since $\Hom_{R\ltimes M}(\textbf{T}(P^{i}),\textbf{Z}(Q))\cong\Hom_{R}(P^{i},Q)$, we have $\Hom_{R\ltimes M}(\Delta,\textbf{Z}(Q))\cong\Hom_{R}(\Xi,Q)$  is exact. Since $\Hom_{R\ltimes M}(\textbf{T}(P^{i}),\textbf{Z}(M\otimes_{R}Q))\cong\Hom_{R}(P^{i},M\otimes_{R}Q)$, $M\otimes_{R}Q\in$ Add$(M)$ and $M$ is a generalized  compatible $R$-$R$-bimodule, we have $\Hom_{R\ltimes M}(\Delta,\textbf{Z}(M\otimes_{R}Q))\cong\Hom_{R}(\Xi,M\otimes_{R}Q)$ is exact. So $\Hom_{R\ltimes M}(\Delta,\textbf{T}(Q))$ is exact by \cite[Theorem 1.4.7]{EJ}.

It follows that  $(X,\alpha)$ is a Gorenstein projective left $R\ltimes M$-module.
\end{proof}
Specializing $M=R$ in Theorem \ref{thm: 3.3}, we have
\begin{cor} \label{cor: 3.6}Let $(X,\alpha)$ be a left $R\ltimes R$-module. If the sequence $X\stackrel{\alpha}\rightarrow X\stackrel{\alpha}\rightarrow X$ is exact and $\coker(\alpha)$ is a Gorenstein projective left $R$-module, then $(X,\alpha)$ is a Gorenstein projective left $R\ltimes R$-module.
\end{cor}
\begin{thm} \label{thm: 3.4}Let $\textbf{Z}(R)$ be a generalized  compatible $R\ltimes M$-$R\ltimes M$-bimodule and $(X,\alpha)$ be a Gorenstein projective left $R\ltimes M$-module. Then the sequence $M\otimes_R M\otimes_R X\stackrel{M\otimes\alpha}\rightarrow M\otimes_R X\stackrel{\alpha}\rightarrow X$ is exact and $\coker(\alpha)$ is a Gorenstein projective left $R$-module.
\end{thm}
\begin{proof}There is a complete projective resolution in $R\ltimes M$-Mod $$\Delta: \cdots\rightarrow \textbf{T}(P^{-1})\rightarrow
\textbf{T}(P^{0})\stackrel{g^{0}}\rightarrow \textbf{T}(P^{1})\rightarrow \textbf{T}(P^{2})\rightarrow \cdots$$
with $(X,\alpha)\cong\ker(g^{0})$.

Since $\textbf{Z}(R)$ is a generalized  compatible $R\ltimes M$-$R\ltimes M$-bimodule, $\textbf{Z}(R)\otimes_{R\ltimes M}\Delta$ is exact. By Lemma \ref{lem: 2.2}(1), we have $$\textbf{Z}(R)\otimes_{R\ltimes M} \textbf{T}(P^{i})\cong R\otimes_{R}P^{i}\cong P^{i}.$$ Hence
we get the exact sequence of projective left $R$-modules
$$\textbf{C}(\Delta): \cdots\rightarrow P^{-1}\rightarrow P^{0}\stackrel{\textbf{C}(g^{0})}\rightarrow P^{1}\rightarrow P^{2}\rightarrow \cdots$$ with
$\coker(\alpha)\cong\ker(\textbf{C}(g^{0}))$.

Let $G$ be a projective left $R$-module. Then $\textbf{Z}(G)\cong$ Add$(\textbf{Z}(R))$. Since $\textbf{Z}(R)$ is a  generalized compatible $R\ltimes M$-$R\ltimes M$-bimodule, we have $\Hom_{R}(\textbf{C}(\Delta),G)\cong\Hom_{R\ltimes M}(\Delta, \textbf{Z}(G))$ is exact. So $\coker(\alpha)$ is a Gorenstein projective left $R$-module.

On the other hand, by Lemma \ref{lem: 2.2}(3), there is an exact sequence in Mod-$R\ltimes M$ $$0\rightarrow \textbf{Z}(M)\rightarrow \textbf{T}(R)\rightarrow \textbf{Z}(R)\rightarrow 0,$$ which induces the exact sequence of complexes  $$0\rightarrow \textbf{Z}(M)\otimes_{R\ltimes M}\Delta\rightarrow \textbf{T}(R)\otimes_{R\ltimes M}\Delta\rightarrow \textbf{Z}(R)\otimes_{R\ltimes M}\Delta\rightarrow 0.$$ Since $\textbf{T}(R)\otimes_{R\ltimes M}\Delta$ and $\textbf{Z}(R)\otimes_{R\ltimes M}\Delta$ are exact, $\textbf{Z}(M)\otimes_{R\ltimes M}\Delta$ is exact  by \cite[Theorem 1.4.7]{EJ}. So $M\otimes_{R}\textbf{C}(\Delta)\cong\textbf{Z}(M)\otimes_{R\ltimes M}\Delta$ is exact by Lemma \ref{lem: 2.2}(1). Let $\iota: \coker(\alpha)\rightarrow P^{0}$ be the injection. Then $M\otimes \iota: M\otimes_{R}\coker(\alpha)\rightarrow M\otimes_{R}P^{0}$ is a monomorphism.

Let $\rho: X\rightarrow \coker(\alpha)$ be the canonical epimorphism. By Lemma \ref{lem: 2.2}(5), there is $\delta:  M\otimes_{R}\coker(\alpha)\rightarrow X$ such that $\delta(M\otimes\rho)=\alpha$. Let $\lambda: X\rightarrow P^{0}\oplus (M\otimes_{R}P^{0})$  and $\varphi^{0}: M\otimes_{R}P^{0}\rightarrow P^{0}\oplus (M\otimes_{R}P^{0})$ be the injections.
By \cite[p.57]{FGR}, we get the following commutative diagram in $R$-Mod: $$\xymatrix{M\otimes_{R}\coker(\alpha)\ar[d]_{\delta}\ar[rr]^{M\otimes \iota}&& M\otimes_{R}P^{0}\ar[d]^{\varphi^{0}}\\
X\ar[rr]^{\lambda}&&P^{0}\oplus (M\otimes_{R}P^{0}).}$$ Then $\delta$ is a monomorphism. Since the sequence $$M\otimes_R M\otimes_R X\stackrel{M\otimes\alpha}\rightarrow M\otimes_R X\stackrel{M\otimes\rho}\rightarrow M\otimes_{R}\coker(\alpha)\rightarrow 0$$ is exact, the sequence $$M\otimes_R M\otimes_R X\stackrel{M\otimes\alpha}\rightarrow M\otimes_R X\stackrel{\alpha}\rightarrow X$$ is also exact.
\end{proof}
Combining Theorem \ref{thm: 3.3} with Theorem \ref{thm: 3.4}, we have
\begin{cor} \label{cor: 3.5}Let $M$ be a  generalized compatible $R$-$R$-bimodule and $\textbf{Z}(R)$ be a generalized  compatible $R\ltimes M$-$R\ltimes M$-bimodule. \begin{enumerate}\item $(X,\alpha)$ is a Gorenstein projective left $R\ltimes M$-module if and only if the sequence $M\otimes_R M\otimes_R X\stackrel{M\otimes\alpha}\rightarrow M\otimes_R X\stackrel{\alpha}\rightarrow X$ is exact and $\coker(\alpha)$ is a Gorenstein projective left $R$-module.\item $\textbf{T}(X)$ is a Gorenstein projective left $R\ltimes M$-module if and only if  $X$ is a Gorenstein projective left $R$-module.\item $\textbf{Z}(X)$ is a Gorenstein projective left $R\ltimes M$-module if and only if $M\otimes_R X=0$ and $X$ is a Gorenstein projective left $R$-module.
\end{enumerate}
\end{cor}
At the end of this section, we construct an example of trivial ring extensions $R\ltimes M$ such that $\textbf{Z}(R)$ is a generalized  compatible $R\ltimes M$-$R\ltimes M$-bimodule.
\begin{exa} \label{exa: 3.7}{\rm Let $R$ be a commutative perfect ring and $M$ be a projective $R$-module. Then we obtain the formal triangular matrix ring $\Gamma=\biggl(\begin{matrix} R&0\\M&R \end{matrix}\biggr)$  with the usual matrix addition and multiplication. It is well known that the ring $\Gamma=\biggl(\begin{matrix} R&0\\M&R \end{matrix}\biggr)$  is exactly the trivial ring extension $(R\times R) \ltimes M$, where $M$ attains an $R\times R$-$R\times R$-bimodule structure through the obvious ring homomorphisms such that $M\otimes_{R\times R}M=0$.
Note that $(R\times R)\ltimes M$ is a perfect ring by \cite[Proposition 1.15]{FGR}. So $pd(\textbf{Z}(R\times R))<\infty$  by \cite[Corollary 4.7]{FGR}. Thus $\textbf{Z}(R\times R)$ is a generalized  compatible $(R\times R)\ltimes M$-$(R\times R)\ltimes M$-bimodule by Example \ref{exa: 3.2}.}
\end{exa}
\bigskip
\section {\bf Gorenstein injective and flat modules over trivial ring extensions}
\bigskip
Recall that an exact sequence of injective left $R$-modules $$\Xi: \cdots\rightarrow
E^{-1}\stackrel{f^{-1}}\rightarrow
E^{0}\stackrel{f^{0}}\rightarrow
E^{1}\stackrel{f^{1}}\rightarrow
E^{2}\rightarrow \cdots$$ is a \emph{complete injective resolution} in $R$-Mod if $\Hom_{R}(Q,\Xi)$ is also exact for any injective left $R$-module $Q$.

A left $R$-module $X$ is called \emph{Gorenstein injective} \cite{EJ1} if there is a complete injective resolution $\Xi: \cdots\rightarrow
E^{-1}\stackrel{f^{-1}}\rightarrow
E^{0}\stackrel{f^{0}}\rightarrow
E^{1}\stackrel{f^{1}}\rightarrow
E^{2}\rightarrow \cdots$ in $R$-Mod  such that $X\cong\ker(f^{0})$.
\begin{df} \label{df: 4.1} {\rm $_{A}N_{B}$ is called a \emph{generalized cocompatible bimodule} if the following two conditions hold:\\
(1) If $\Xi$ is a complete injective resolution  in $A$-Mod, then $\Hom_{A}(N,\Xi)$ is exact.\\
(2) If $\Delta$ is a complete injective resolution in $B$-Mod, then $\Hom_{B}$(Prod$(N^{+}),\Delta)$ is exact.}
\end{df}
\begin{exa} \label{exa: 4.2}{\rm (1) If $fd(N_{B})<\infty$ and $pd(_{A}N)<\infty$, then $_{A}N_{B}$ is a generalized cocompatible bimodule by \cite[Lemma 2.5]{EIT} and \cite[Theorem 2.1]{F}.

(2) If $fd(N_{B})<\infty$ and $id(_{A}N)<\infty$, then $_{A}N_{B}$ is a generalized cocompatible bimodule by \cite[Theorem 2.1]{F}. }
\end{exa}
\begin{thm} \label{thm: 4.3}Let $M$ be a generalized  cocompatible $R$-$R$-bimodule and $[Y,\beta]$ be a left $R\ltimes M$-module such that the sequence $Y\stackrel{\beta}\rightarrow \Hom_R(M, Y)\stackrel{\beta_{*}}\rightarrow\Hom_R(M,\Hom_R(M, Y))$ is exact and $\ker(\beta)$ is a Gorenstein injective left $R$-module. Then $[Y,\beta]$  is a Gorenstein injective left $R\ltimes M$-module.
\end{thm}
\begin{proof}There exists a complete injective resolution in $R$-Mod $$\Xi: \cdots\rightarrow
E^{-1}\stackrel{f^{-1}}\rightarrow
E^{0}\stackrel{f^{0}}\rightarrow
E^{1}\stackrel{f^{1}}\rightarrow
E^{2}\rightarrow \cdots$$  with $\ker(\beta)\cong\ker(f^{0})$. Since $M$ is a generalized  cocompatible $R$-$R$-bimodule, we get the exact sequence of left $R$-modules
$$\Hom_{R}(M,\Xi): \cdots\rightarrow\Hom_{R}(M,E^{-1})\stackrel{(f^{-1})_{*}}\rightarrow
\Hom_{R}(M,E^{0})\stackrel{(f^{0})_{*}}\rightarrow
\Hom_{R}(M,E^{1})\rightarrow \cdots$$ with
$\Hom_{R}(M,\ker(\beta))\cong\ker((f^{0})_{*})$. Since the sequences $$0\rightarrow \ker(\beta)\stackrel{\iota}\rightarrow Y\stackrel{\beta}\rightarrow \Hom_R(M, Y)\stackrel{\beta_{*}}\rightarrow \Hom_R(M,\Hom_R(M, Y))$$ and $$0\rightarrow \Hom_{R}(M,\ker(\beta))\stackrel{\iota_{_{*}}}\rightarrow \Hom_R(M, Y)\stackrel{\beta_{*}}\rightarrow \Hom_R(M,\Hom_R(M, Y))$$ are  exact, we get the exact sequence $$0\rightarrow \ker(\beta)\stackrel{\iota}\rightarrow Y\stackrel{\gamma}\rightarrow \Hom_R(M,\ker(\beta))\rightarrow 0$$ such that $\beta=\iota_{*}\gamma$ by Lemma \ref{lem: 2.2}(6).

There exists a split monomorphism $E^{i}\rightarrow \prod R^{+}$ for any $i$, which induces the split monomorphism $$\Hom_{R}(M,E^{i})\rightarrow \Hom_{R}(M,\prod R^{+})\cong \prod M^{+}.$$ So $\Hom_{R}(M,E^{i})\in$ Prod$(M^{+})$.  Since $M$ is a generalized cocompatible $R$-$R$-bimodule, $\Hom_{R}(\Hom_{R}(M,E^{i}),\Xi)$ is exact. Thus $\Ext^{1}_{R}(\Hom_{R}(M, E^{i}),\ker(f^{i+1}))=0$.

Let $\pi: E^{-1}\rightarrow \ker(\beta)$ be the canonical epimorphism and $\tau:\ker(\beta)\rightarrow E^{0}$ be the injection such that $\tau\pi=f^{-1}$. There is $\psi: \Hom_{R}(M,E^{-1})\rightarrow Y$ such that $\gamma\psi=\pi_{*}$  and there is $\phi: Y\rightarrow E^{0}$ such that $\phi\iota=\tau$.

Define $\xi: \Hom_{R}(M,E^{-1})\oplus E^{-1}\rightarrow Y$ by $$\xi(x,y)=\psi(x)+\iota\pi(y)$$ and define $\lambda: Y\rightarrow\Hom_{R}(M,E^{0})\oplus E^{0}$  by $$\lambda(x)=(\tau\gamma(x),\phi(x)).$$ It is easy to check that $\lambda$ is a monomorphism and $\xi$ is an epimorphism.  Then by the generalized Horseshoe Lemma \cite[Lemma 1.6]{Z}, we get the following commutative diagrams:
$$\xymatrix{&\vdots\ar[d]
&&\vdots \ar@{.>}[d] &\vdots \ar[d] &\\
0\ar[r] &E^{-2}\ar[d]^{f^{-2}}\ar[rr]&&\Hom_{R}(M,E^{-2})\oplus E^{-2}\ar@{.>}[d]^{\partial^{-2}}
\ar[r]&\Hom_{R}(M,E^{-2})\ar[r]\ar[d]^{(f^{-2})_{*}}&0\\0
\ar[r]&E^{-1}\ar[d]^{\pi}\ar[rr]&&\Hom_{R}(M,E^{-1})\oplus E^{-1}\ar@{.>}[d]^{\xi}
\ar[r]&\Hom_{R}(M,E^{-1})\ar[r]\ar[d]^{\pi_{*}}&0\\
0 \ar[r] &\ker(\beta)\ar[d]\ar[rr]^{\iota}&&Y\ar[d]\ar[r]^{\gamma}&\Hom_R(M, \ker(\beta))\ar[r]\ar[d]&0\\
& 0&&0 &0 &}$$
and
 $$\xymatrix{&0\ar[d]
&&0\ar[d] &0\ar[d] &\\
0 \ar[r] &\ker(\beta)\ar[d]^{\tau}\ar[rr]^{\iota}&&Y\ar@{.>}[d]^{\lambda}
\ar[r]^{\gamma}&\Hom_R(M, \ker(\beta))\ar[r]\ar[d]^{\tau_{*}}&0\\0
\ar[r] &E^{0}\ar[d]^{f^{0}}\ar[rr]&&\Hom_{R}(M,E^{0})\oplus E^{0}\ar@{.>}[d]^{\partial^{0}}
\ar[r]&\Hom_{R}(M,E^{0})\ar[r]\ar[d]^{(f^{0})_{*}}&0\\0
\ar[r]&E^{1}\ar[d]^{f^{1}}\ar[rr]&&\Hom_{R}(M,E^{1})\oplus E^{1}\ar@{.>}[d]^{\partial^{1}}
\ar[r]&\Hom_{R}(M,E^{1})\ar[r]\ar[d]^{(f^{1})_{*}}&0\\
0\ar[r]&E^{2}\ar[d]\ar[rr]&&\Hom_{R}(M,E^{2})\oplus E^{2}
\ar@{.>}[d]\ar[r]&\Hom_{R}(M,E^{2})\ar[d]\ar[r]&0\\
&\vdots &&\vdots &\vdots &}$$ with exact rows and columns and $\partial^{i}=\biggl(\begin{matrix} (f^{i})_{*}&0\\\tau^{i}&f^{i}\end{matrix}\biggr)$.

Let $\partial^{-1}=\lambda\xi$. Then we get the exact sequence of  left $R$-modules
$$\cdots\rightarrow \Hom_{R}(M,E^{-1})\oplus E^{-1}\stackrel{\partial^{-1}}\rightarrow \Hom_{R}(M,E^{0})\oplus E^{0}\stackrel{\partial^{0}}\rightarrow \Hom_{R}(M,E^{1})\oplus E^{1}\rightarrow\cdots$$  with $Y\cong\ker(\partial^{0})$. The following commutative diagram
$${\Small \xymatrix{\cdots\ar[r]&\Hom_{R}(M,E^{-2})\oplus E^{-2}\ar[d]_{\vartheta_{-2}}\ar[r]^{\partial^{-2}}&\Hom_{R}(M,E^{-1})\oplus E^{-1}\ar[d]^{\vartheta_{-1}}\ar[r]^{\xi}&Y\ar[d]^{\beta}\ar[r]&0\\
\cdots\ar[r]&\Hom_{R}(M,\Hom_{R}(M,E^{-2})\oplus E^{-2})\ar[r]^{(\partial^{-2})_{*}}&\Hom_{R}(M,\Hom_{R}(M,E^{-1})\oplus E^{-1})\ar[r]^{\xi_{*}}&\Hom_{R}(M,Y)}}$$ implies that the sequence $\cdots\rightarrow\textbf{H}(E^{-2})\stackrel{\partial^{-2}}\rightarrow \textbf{H}(E^{-1})\stackrel{\xi}\rightarrow[Y,\beta]\rightarrow 0$ is exact.

The following commutative diagram
$${\Small \xymatrix{0\ar[r]&Y\ar[d]_{\beta}\ar[rr]^{\lambda}&&\Hom_{R}(M,E^{0})\oplus E^{0}\ar[d]^{\vartheta_{0}}\ar[r]^{\partial^{0}}&\Hom_{R}(M,E^{1})\oplus E^{1}\ar[d]^{\vartheta_{1}}\ar[r]&\cdots\\
&\Hom_{R}(M,Y)\ar[rr]^{\lambda_{*}}&&\Hom_{R}(M,\Hom_{R}(M,E^{0})\oplus E^{0})\ar[r]^{(\partial^{0})_{*}}&\Hom_{R}(M,\Hom_{R}(M,E^{1})\oplus E^{1})\ar[r]&\cdots}}$$ implies that the sequence $0\rightarrow [Y,\beta]\stackrel{\lambda}\rightarrow \textbf{H}(E^{0})\stackrel{\partial^{0}}\rightarrow \textbf{H}(E^{1})\rightarrow\cdots$ is exact.

By Lemma \ref{lem: 2.1}(2), each $\textbf{H}(E^{i})$ is injective. So  there is an exact sequence of injective left $R\ltimes M$-modules $$\Delta: \cdots\rightarrow \textbf{H}(E^{-1})\stackrel{\partial^{-1}}\rightarrow
\textbf{H}(E^{0})\stackrel{\partial^{0}}\rightarrow \textbf{H}(E^{1})\stackrel{\partial^{1}}\rightarrow \textbf{H}(E^{2})\rightarrow \cdots$$ with $[Y,\beta]\cong\ker(\partial^{0})$.

By Lemma \ref{lem: 2.1}(2), any injective  left $R\ltimes M$-module is isomorphic to $\textbf{H}(Q)$ with $Q$  an injective  left $R$-module. By Lemma \ref{lem: 2.2}(4), there is an exact sequence in $R\ltimes M$-Mod $$0\rightarrow \textbf{Z}(Q)\rightarrow \textbf{H}(Q)\rightarrow \textbf{Z}(\Hom_{R}(M,Q))\rightarrow 0,$$ which induces the exact sequence  of complexes $$0\rightarrow \Hom_{R\ltimes M}(\textbf{Z}(\Hom_{R}(M,Q)),\Delta)\rightarrow \Hom_{R\ltimes M}(\textbf{H}(Q),\Delta)\rightarrow \Hom_{R\ltimes M}(\textbf{Z}(Q),\Delta)\rightarrow 0.$$
Note that $\Hom_{R\ltimes M}(\textbf{Z}(Q),\textbf{H}(E^{i}))\cong\Hom_{R}(Q,E^{i})$, so the complex $\Hom_{R\ltimes M}(\textbf{Z}(Q),\Delta)\cong\Hom_{R}(Q,\Xi)$  is exact. Since $\Hom_{R\ltimes M}(\textbf{Z}(\Hom_{R}(M,Q)),\textbf{H}(E^{i}))\cong\Hom_{R}(\Hom_{R}(M,Q),E^{i})$, $\Hom_{R}(M,Q)\in$ Prod$(M^{+})$ and $M$ is a generalized cocompatible $R$-$R$-bimodule, the complex $\Hom_{R\ltimes M}(\textbf{Z}(\Hom_{R}(M,Q)),\Delta)\cong\Hom_{R}(\Hom_{R}(M,Q),\Xi)$ is exact. Therefore $\Hom_{R\ltimes M}(\textbf{H}(Q),\Delta)$ is exact  by \cite[Theorem 1.4.7]{EJ}.

It follows that $[Y,\beta]$ is a Gorenstein injective left $R\ltimes M$-module.
\end{proof}
Specializing $M=R$ in Theorem \ref{thm: 4.3}, we have
\begin{cor}\label{cor: 4.6}Let $[Y,\beta]$  be a left $R\ltimes R$-module. If the sequence $Y\stackrel{\beta}\rightarrow Y\stackrel{\beta}\rightarrow Y$ is exact and $\ker(\beta)$ is a Gorenstein injective left $R$-module, then $[Y,\beta]$  is a Gorenstein injective left $R\ltimes R$-module.
\end{cor}
\begin{thm}\label{thm: 4.4}Let $\textbf{Z}(R)$ be a generalized cocompatible $R\ltimes M$-$R\ltimes M$-bimodule and $[Y,\beta]$ be a Gorenstein injective left $R\ltimes M$-module. Then the sequence $Y\stackrel{\beta}\rightarrow\Hom_R(M, Y)\stackrel{\beta_{*}}\rightarrow\Hom_R(M,\Hom_R(M, Y))$ is exact and $\ker(\beta)$ is a Gorenstein injective left $R$-module.
\end{thm}
\begin{proof}There is a complete injective resolution in $R\ltimes M$-Mod $$\Delta: \cdots\rightarrow \textbf{H}(E^{-1})\rightarrow
\textbf{H}(E^{0})\stackrel{\partial^{0}}\rightarrow \textbf{H}(E^{1})\rightarrow \textbf{H}(E^{2})\rightarrow \cdots$$
with $[Y,\beta]\cong\ker(\partial^{0})$. Since $\textbf{Z}(R)$ is a generalized cocompatible $R\ltimes M$-$R\ltimes M$-bimodule, $\Hom_{R\ltimes M}(\textbf{Z}(R),\Delta)$ is exact. Note that $\Hom_{R\ltimes M}(\textbf{Z}(R),[Y,\beta])\cong\Hom_{R}(R,\ker(\beta))\cong \ker(\beta)$ by Lemma \ref{lem: 2.2}(2). So we get the exact sequence of injective left $R$-modules
$$\textbf{K}(\Delta): \cdots\rightarrow E^{-1}\rightarrow E^{0}\stackrel{\textbf{K}(\partial^{0})}\rightarrow E^{1}\rightarrow E^{2}\rightarrow \cdots$$  with $\ker(\beta)\cong\ker(\textbf{K}(\partial^{0}))$.

Let $G$ be an injective left $R$-module. Then there is a split monomorphism $G\rightarrow \prod R^{+}$, which induces the split monomorphism $$\textbf{Z}(G)\rightarrow\textbf{Z}(\prod R^{+})\cong \prod\textbf{Z}(R^{+}).$$  So $\textbf{Z}(G)\in$ Prod$(\textbf{Z}(R)^{+})$.  Since $\textbf{Z}(R)$ is a generalized cocompatible $R\ltimes M$-$R\ltimes M$-bimodule, $\Hom_{R\ltimes M}(\textbf{Z}(G),\Delta)$ is exact. By Lemma \ref{lem: 2.2}(2), $\Hom_{R\ltimes M}(\textbf{Z}(G),\textbf{H}(E^{i}))\cong \Hom_{R}(G,E^{i})$. Hence $\Hom_{R}(G,\textbf{K}(\Delta))\cong\Hom_{R\ltimes M}(\textbf{Z}(G),\Delta)$ is exact. So $\ker(\beta)$ is a Gorenstein injective left $R$-module.

By Lemma \ref{lem: 2.2}(3), there is an exact sequence in $R\ltimes M$-Mod  $$0\rightarrow \textbf{Z}(M)\rightarrow \textbf{T}(R)\rightarrow \textbf{Z}(R)\rightarrow 0,$$ which induces the exact sequence of complexes  $$0\rightarrow \Hom_{R\ltimes M}(\textbf{Z}(R),\Delta)\rightarrow \Hom_{R\ltimes M}(\textbf{T}(R),\Delta)\rightarrow \Hom_{R\ltimes M}(\textbf{Z}(M),\Delta)\rightarrow 0.$$ Since both $\Hom_{R\ltimes M}(\textbf{Z}(R),\Delta)$  and $\Hom_{R\ltimes M}(\textbf{T}(R),\Delta)$ are exact, $\Hom_{R\ltimes M}(\textbf{Z}(M),\Delta)$ is exact  by \cite[Theorem 1.4.7]{EJ}. So $\Hom_{R}(M,\textbf{K}(\Delta))\cong\Hom_{R\ltimes M}(\textbf{Z}(M),\Delta)$ is exact.
Let $\tau: E^{-1}\rightarrow\ker(\beta)$  be the canonical epimorphism. Then $\tau_{*}: \Hom_{R}(M,E^{-1})\rightarrow\Hom_{R}(M,\ker(\beta))$ is an epimorphism.
Let $\iota: \ker(\beta)\rightarrow Y$ be the inclusion. By Lemma \ref{lem: 2.2}(6), there is $\gamma: Y\rightarrow \Hom_{R}(M,\ker(\beta))$ such that $\iota_{*}\gamma=\beta$. Let $\rho: \Hom_{R}(M,E^{-1})\oplus E^{-1}\rightarrow Y$  be the canonical epimorphism and $\varphi^{-1}:\Hom_{R}(M,E^{-1})\oplus E^{-1}\rightarrow\Hom_{R}(M,E^{-1})$ be the projection.
Consider the following commutative diagram: $$\xymatrix{\Hom_{R}(M,E^{-1})\oplus E^{-1}\ar[d]^{\varphi^{-1}}\ar[rr]^{\rho}&&Y\ar[d]^{\gamma}\\
\Hom_{R}(M,E^{-1})\ar[rr]^{\tau_{*}}&&\Hom_{R}(M,\ker(\beta)).}$$
Since $\tau_{*}$  and $\varphi^{-1}$ are epimorphisms, $\gamma$ is an epimorphism. Hence the sequence  $Y\stackrel{\beta}\rightarrow \Hom_R(M, Y)\stackrel{\beta_{*}}\rightarrow \Hom_R(M,\Hom_R(M, Y))$ is exact.
\end{proof}
Combining Theorem \ref{thm: 4.3} with Theorem \ref{thm: 4.4}, we have
\begin{cor} \label{cor: 4.5}Let $M$ be a  generalized cocompatible $R$-$R$-bimodule and $\textbf{Z}(R)$ be a generalized  cocompatible $R\ltimes M$-$R\ltimes M$-bimodule.\begin{enumerate}\item  $[Y,\beta]$  is a Gorenstein injective left $R\ltimes M$-module if and only if the sequence  $Y\stackrel{\beta}\rightarrow \Hom_R(M, Y)\stackrel{\beta_{*}}\rightarrow \Hom_R(M,\Hom_R(M, Y))$ is   exact and $\ker(\beta)$ is a Gorenstein injective left $R$-module.\item $\textbf{H}(Y)$  is a Gorenstein injective left $R\ltimes M$-module if and only if  $Y$ is a Gorenstein injective left $R$-module.\item $\textbf{Z}(Y)$  is a Gorenstein injective left $R\ltimes M$-module if and only if $\Hom_R(M, Y)=0$ and $Y$ is a Gorenstein injective left $R$-module.
\end{enumerate}
\end{cor}
We note that Example  \ref{exa: 3.7} is also an example of trivial ring extensions $R\ltimes M$ such that $\textbf{Z}(R)$ is a generalized  cocompatible $R\ltimes M$-$R\ltimes M$-bimodule by Example \ref{exa: 4.2}.

At the end of this section, we describe Gorenstein flat modules over trivial ring extensions.

Recall that a right $R$-module $X$ is \emph{Gorenstein flat} \cite{EJT} if there is an exact sequence $\cdots\rightarrow F^{-1}\rightarrow
F^{0}\stackrel{f^{0}}\rightarrow F^{1}\rightarrow F^{2}\rightarrow \cdots$  of flat right $R$-modules with $X\cong\ker(f^{0})$, which remains exact after applying $-\otimes_{R}E$ for any injective left $R$-module $E$.

Recall that $R$ is  a \emph{left coherent ring} \cite{L} if every
finitely generated  left ideal of $R$ is finitely presented.
\begin{thm} \label{thm: 4.7}Let  $R\ltimes M$ be a left coherent ring and $(X,\alpha)$ be a right $R\ltimes M$-module.
 \begin{enumerate}\item If $M$ is a generalized   cocompatible $R$-$R$-bimodule, the sequence $X\otimes_R M\otimes_R M \stackrel{\alpha\otimes M}\rightarrow X\otimes_R M \stackrel{\alpha}\rightarrow X$ is exact and $\coker(\alpha)$ is a Gorenstein flat right $R$-module, then $(X,\alpha)$ is a Gorenstein flat right $R\ltimes M$-module.\item If $\textbf{Z}(R)$ is a generalized cocompatible $R\ltimes M$-$R\ltimes M$-bimodule  and $(X,\alpha)$ is a  Gorenstein flat right $R\ltimes M$-module, then the sequence $X\otimes_R M\otimes_R M \stackrel{\alpha\otimes M}\rightarrow X\otimes_R M \stackrel{\alpha}\rightarrow X$  is exact and $\coker(\alpha)$ is a Gorenstein flat right $R$-module.
 \end{enumerate}
\end{thm}
\begin{proof}We note that $R$ is also a left coherent ring by \cite[Theorem 2.2]{FGR}.

(1) Note that the sequence $X^{+}\stackrel{\alpha^{+}}\rightarrow\Hom_R(M, X^{+})\stackrel{(\alpha^{+})_{*}}\rightarrow \Hom_R(M,\Hom_R(M, X^{+}))$ is   exact and $\ker(\alpha^{+})\cong (\coker(\alpha))^{+}$ is a Gorenstein injective left $R$-module by \cite[Theorem 3.6]{H}. Thus $(X,\alpha)^{+}\cong[X^{+},\alpha^{+}]$  is a Gorenstein injective left $R\ltimes M$-module by Theorem \ref{thm: 4.3}. So $(X,\alpha)$ is a Gorenstein flat right $R\ltimes M$-module by \cite[Theorem 3.6]{H}.

(2) Note that $[X^{+},\alpha^{+}]\cong(X,\alpha)^{+}$  is a Gorenstein injective left $R\ltimes M$-module by \cite[Theorem 3.6]{H}. So the sequence $X^{+}\stackrel{\alpha^{+}}\rightarrow\Hom_R(M, X^{+})\stackrel{(\alpha^{+})_{*}}\rightarrow \Hom_R(M,\Hom_R(M, X^{+}))$  is exact and $(\coker(\alpha))^{+}\cong\ker(\alpha^{+})$ is a Gorenstein injective left $R$-module by Theorem \ref{thm: 4.4}. Hence the sequence $X\otimes_R M\otimes_R M \stackrel{\alpha\otimes M}\rightarrow X\otimes_R M \stackrel{\alpha}\rightarrow X$  is exact and $\coker(\alpha)$ is a Gorenstein flat right $R$-module by \cite[Theorem 3.6]{H}.
\end{proof}
\begin{cor}\label{cor: 4.8}Let $R\ltimes M$ be a left coherent ring, $M$ be a generalized cocompatible $R$-$R$-bimodule and $\textbf{Z}(R)$ be a generalized  cocompatible $R\ltimes M$-$R\ltimes M$-bimodule.
\begin{enumerate}\item $(X,\alpha)$ is a Gorenstein flat right $R\ltimes M$-module if and only if the sequence $X\otimes_R M\otimes_R M \stackrel{\alpha\otimes M}\rightarrow X\otimes_R M \stackrel{\alpha}\rightarrow X$ is exact and $\coker(\alpha)$ is a Gorenstein flat right $R$-module.\item $\textbf{T}(X)$ is a Gorenstein flat right $R\ltimes M$-module if and only if  $X$ is a Gorenstein flat right $R$-module.\item $\textbf{Z}(X)$ is a Gorenstein flat right $R\ltimes M$-module if and only if $X\otimes_R M =0$ and $X$ is a Gorenstein flat right $R$-module.
\end{enumerate}
\end{cor}
Specializing $M=R$ in Theorem \ref{thm: 4.7}, we have
\begin{cor}\label{cor: 4.9}Let  $R$ be a left coherent ring and $(X,\alpha)$   be a right $R\ltimes R$-module. If the sequence $X\stackrel{\alpha}\rightarrow X\stackrel{\alpha}\rightarrow X$ is exact and $\coker(\alpha)$ is a Gorenstein flat right $R$-module, then $(X,\alpha)$  is a Gorenstein flat right $R\ltimes R$-module.
\end{cor}
\bigskip
\section {\bf Applications}
\bigskip
 Morita context rings have been studied explicitly in numerous papers and books \cite{FGR,GP,G1,GP1,KT,M1,YY}. For example, Green and Psaroudakis gave bounds for the global dimension of a Morita context  ring  which is an Artin algebra and has  zero bimodule homomorphisms \cite{GP1}. Yan and Yao studied pure projectivity and $FP$-injectivity of modules over a Morita context ring with  zero bimodule homomorphisms \cite{YY}. Gao and Psaroudakis investigated how to construct Gorenstein projective modules over a Morita context ring which is an Artin algebra and has  zero bimodule homomorphisms \cite{GP}.

In this section, we apply the foregoing results to Morita context rings with  zero bimodule homomorphisms since this kind of rings is one special case of trivial ring extensions.

Let $A$ and $B$ be rings, $_{B}U_{A}$ and $_{A}V_{B}$ be bimoduls, $\phi: U\otimes _{A}V\rightarrow B$ and $\psi: V\otimes_{B}U\rightarrow A$ be bimodule homomorphisms. Then $\biggl(\begin{matrix} A&_{A}V_{B}\\_{B}U_{A}&B \end{matrix}\biggr)_{(\phi,\psi)}$ becomes a ring with usual matrix addition and multiplication given by
$$\biggl(\begin{matrix} a_{1}&v_{1}\\u_{1}&b_{1}\end{matrix}\biggr)\biggl(\begin{matrix} a_{2}&v_{2}\\u_{2}&b_{2} \end{matrix}\biggr)=\biggl(\begin{matrix} a_{1}a_{2}+\psi(v_{1}\otimes u_{2})&a_{1}v_{2}+v_{1}b_{2}\\u_{1}a_{2}+b_{1}u_{2}&b_{1}b_{2}+\phi(u_{1}\otimes v_{2})\end{matrix}\biggr).$$ $\biggl(\begin{matrix} A&_{A}V_{B}\\_{B}U_{A}&B \end{matrix}\biggr)_{(\phi,\psi)}$  is called  a \emph{Morita context ring} \cite{M1} or \emph{formal matrix ring} \cite{KT}. In particular, if $\phi=0,\psi=0$, then $\biggl(\begin{matrix} A&_{A}V_{B}\\_{B}U_{A}&B \end{matrix}\biggr)_{(0,0)}$ is called  a \emph{Morita context ring with  zero bimodule homomorphisms}, which is a generalization of the formal triangular matrix ring $\biggl(\begin{matrix} A&0\\_{B}U_{A}&B \end{matrix}\biggr)$.

Let $\Lambda=\biggl(\begin{matrix} A&_{A}V_{B}\\_{B}U_{A}&B \end{matrix}\biggr)_{(0,0)}$. Green \cite{G1} proved that the category $\Lambda$-Mod is equivalent to the category
$\Omega$ whose objects are tuples $(X,Y,f,g)$, where $X\in A$-Mod, $Y\in B$-Mod, $f\in\Hom_{B}(U\otimes_{A} X, Y)$ and $g\in\Hom_{A}(V\otimes_{B} Y, X)$ such that $g(V\otimes f)=0, f(U\otimes g)=0$ and whose
morphisms from $(X_{1},Y_{1},f_{1},g_{1})$ to $(X_{2},Y_{2},f_{2},g_{2})$ are pairs $(\alpha,\beta)$  such that $\alpha\in\Hom_{A}(X_{1}, X_{2}), \beta\in\Hom_{B}(Y_{1}, Y_{2})$ and
the following diagrams are commutative.
$$\xymatrix{U\otimes_{A} X_{1}\ar[d]_{f_{1}}\ar[r]^{U\otimes \alpha}&U\otimes_{A} X_{2}\ar[d]^{f_{2}}\\
Y_{1}\ar[r]^{\beta}&Y_{2}}\xymatrix{V\otimes_{B} Y_{1}\ar[d]_{g_{1}}\ar[r]^{V\otimes \beta}&V\otimes_{B} Y_{2}\ar[d]^{g_{2}}\\
X_{1}\ar[r]^{\alpha}&X_{2}}$$ In view of the well-known adjointness relation, the category $\Lambda$-Mod is also equivalent to the category $\Gamma$ whose objects are tuples $[X,Y,f,g]$, where $X\in A$-Mod, $Y\in B$-Mod, $f\in\Hom_{A}(X,\Hom_{B}(U, Y))$ and $g\in \Hom_{B}(Y,\Hom_{A}(V, X))$ such that $\Hom_{B}(U,g)f=0, \Hom_{A}(U,f)g=0$ and whose morphisms from $[X_{1},Y_{1},f_{1},g_{1}]$ to $[X_{2},Y_{2},f_{2},g_{2}]$ are pairs $[\alpha,\beta]$  such that
$\alpha\in\Hom_{A}(X_{1}, X_{2}), \beta\in\Hom_{B}(Y_{1}, Y_{2})$ and $f_{2}\alpha=\Hom_{B}(U,\beta)f_{1}$, $g_{2}\beta=\Hom_{A}(V,\alpha)g_{1}$.
In what follows, we will identify $\Lambda$-Mod  with $\Omega$ and $\Gamma$.

Note that a sequence $0\rightarrow(X',Y',f',g')\rightarrow(X,Y,f,g)\rightarrow(X'',Y'',f'',g'')\rightarrow 0$ of left
$\Lambda$-modules is exact if and only if both sequences $0\rightarrow
X'\rightarrow X\rightarrow X''\rightarrow 0$ of left
$A$-modules and $0\rightarrow Y'\rightarrow Y\rightarrow Y''\rightarrow
0$ of left $B$-modules are exact.

Analogously, the category Mod-$\Lambda$ of right $\Lambda$-modules is equivalent to the category whose objects are  tuples $(W,Q,f,g)$, where $W\in$ Mod-$A$, $Q\in$ Mod-$B$, $f\in\Hom_{A}(Q\otimes_{B}U,W)$, $g\in\Hom_{B}( W\otimes_{A}V,Q)$ such that $g(f\otimes V)=0, f(g\otimes U)=0$  and whose morphisms from $(W_{1},Q_{1},f_{1},g_{1})$ to $(W_{2},Q_{2},f_{2},g_{2})$ are
pairs $(\alpha,\beta)$  such that $\alpha\in\Hom_{A}(W_{1}, W_{2}),
\beta\in\Hom_{B}(Q_{1}, Q_{2})$  and $g_{2}(\alpha\otimes V)=\beta g_{1}$, $f_{2}(\beta\otimes U)=\alpha f_{1}$.

 Note that $U\oplus V$ attains the $A\times B$-$A\times B$-bimodule structure through the ring homomorphisms $A\times B \rightarrow A$ and $A\times B \rightarrow B$. It is well known that a left $A\times B$-module is an order pair $(X,Y)$ with $X\in A$-Mod and $Y\in B$-Mod. Similarly, a right $A\times B$-module is an order pair $(W_{1},W_{2})$ with $W_{1}\in$ Mod-$A$ and $W_{2}\in$ Mod-$B$. So $(U\oplus V)\otimes_{A\times B} (X,Y)\cong (V\otimes_{B}Y,U\otimes_{A}X)$ and $\Hom_{A\times B}(U\oplus V, (X,Y))\cong(\Hom_{B}(U,Y),\Hom_{A}(V,X))$.  It is known that the ring $\Lambda=\biggl(\begin{matrix} A&_{A}V_{B}\\_{B}U_{A}&B \end{matrix}\biggr)_{(0,0)}$ is isomorphic to the trivial ring extension $(A\times B)\ltimes (U\oplus V)$ under the correspondence: $\biggl(\begin{matrix} a&v\\u&b \end{matrix}\biggr)\rightarrow ((a,b),(u,v))$ \cite{FGR, GP1}. Therefore $\Lambda$-Mod is isomorphic to $(A\times B) \ltimes (U\oplus V)$-Mod by the functor $\Theta: \Lambda$-Mod $\rightarrow (A\times B) \ltimes (U\oplus V)$-Mod given by $\Theta(X,Y,f,g)=((X,Y),(g,f))$. Similarly, Mod-$\Lambda$ is isomorphic to Mod-$(A\times B)\ltimes (U\oplus V)$ by the functor $\Upsilon:$ Mod-$\Lambda\rightarrow$ Mod-$(A\times B)\ltimes (U\oplus V)$ given by $\Upsilon(W,Q,f,g)=((W,Q),(f,g))$.
 \begin{lem} \label{lem: 5.1}Let $_{B}U_{A}$ and $_{A}V_{B}$ be  bimodules.
  \begin{enumerate}\item If both $_{B}U_{A}$ and $_{A}V_{B}$ are generalized compatible bimodules, then $U\oplus V$ is a generalized compatible $A\times B$-$A\times B$-bimodule.\item If both  $_{B}U_{A}$ and $_{A}V_{B}$ are generalized cocompatible bimodules, then $U\oplus V$ is a generalized cocompatible $A\times B$-$A\times B$-bimodule.
\end{enumerate}
\end{lem}
\begin{proof}(1) Given a complete projective resolution in $A\times B$-Mod  $$\Xi: \cdots\rightarrow(P^{-1}_{1},P^{-1}_{2})\rightarrow(P^{0}_{1},P^{0}_{2})\rightarrow(P^{1}_{1},
P^{1}_{2})\rightarrow(P^{2}_{1},P^{2}_{2})\rightarrow\cdots,$$
then $$\Xi_{1}: \cdots\rightarrow P^{-1}_{1}\rightarrow P^{0}_{1}\rightarrow P^{1}_{1}\rightarrow P^{2}_{1}\rightarrow\cdots$$ and $$\Xi_{2}: \cdots\rightarrow P^{-1}_{2}\rightarrow P^{0}_{2}\rightarrow P^{1}_{2}\rightarrow P^{2}_{2}\rightarrow\cdots$$ are  complete projective resolutions in $A$-Mod  and $B$-Mod, respectively.

Since the sequences $U\otimes_{A}\Xi_{1}$ and $V\otimes_{B}\Xi_{2}$  are  exact, the sequence $(U\oplus V)\otimes_{A\times B}\Xi\cong (U\otimes_{A}\Xi_{1})\oplus (V\otimes_{B}\Xi_{2})$ is exact. Since the sequences $\Hom_{A}(\Xi_{1}$, Add$(V))$ and $\Hom_{B}(\Xi_{2}$, Add$(U))$ are exact,  the sequence $\Hom_{A\times B}(\Xi$, Add$(U\oplus V))\cong\Hom_{A}(\Xi_{1}$, Add$(V))\oplus\Hom_{B}(\Xi_{2}$, Add$(U))$ is exact.

Thus $U\oplus V$ is a generalized compatible $A\times B$-$A\times B$-bimodule.

(2)  Given a complete injective resolution in $A\times B$-Mod  $$\Delta: \cdots\rightarrow(E^{-1}_{1},E^{-1}_{2})\rightarrow(E^{0}_{1},
E^{0}_{2})\rightarrow(E^{1}_{1},E^{1}_{2})\rightarrow(E^{2}_{1},E^{2}_{2})\rightarrow\cdots,$$
then $$\Delta_{1}: \cdots\rightarrow E^{-1}_{1}\rightarrow E^{0}_{1}\rightarrow E^{1}_{1}\rightarrow E^{2}_{1}\rightarrow\cdots$$ and $$\Delta_{2}: \cdots\rightarrow E^{-1}_{2}\rightarrow E^{0}_{2}\rightarrow E^{1}_{2}\rightarrow E^{2}_{2}\rightarrow\cdots$$ are  complete injective resolutions in $A$-Mod  and $B$-Mod, respectively.

Since the sequences $\Hom_{A}(V,\Delta_{1})$ and $\Hom_{B}(U,\Delta_{2})$ are  exact, the sequence $\Hom_{A\times B}(U\oplus V,\Delta)\cong \Hom_{A}(V,\Delta_{1})\oplus \Hom_{B}(U,\Delta_{2})$ is exact.

Since the sequences $\Hom_{A}$(Prod$(U^{+}),\Delta_{1})$  and  $\Hom_{B}$(Prod$(V^{+}),\Delta_{2})$  are exact, the sequence  $\Hom_{A\times B}$(Prod$((U\oplus V)^{+}),\Delta)\cong \Hom_{A}$(Prod$(U^{+}),\Delta_{1})\oplus\Hom_{B}$(Prod$(V^{+}),\Delta_{2})$ is exact.

Therefore $U\oplus V$ is a generalized cocompatible $A\times B$-$A\times B$-bimodule.
\end{proof}
\begin{thm} \label{thm: 5.2}Let $\Lambda=\biggl(\begin{matrix} A&_{A}V_{B}\\_{B}U_{A}&B \end{matrix}\biggr)_{(0,0)}$ and $(X,Y,f,g)$ be a left $\Lambda$-module.
 \begin{enumerate}\item If both $_{B}U_{A}$ and $_{A}V_{B}$ are generalized compatible bimodules, the sequences $V\otimes_B U\otimes_A X\stackrel{V\otimes f}\rightarrow V\otimes_B Y\stackrel{g}\rightarrow X$ and $U\otimes_A V\otimes_B Y\stackrel{U\otimes g}\rightarrow U\otimes_A X\stackrel{f}\rightarrow Y$ are exact, $\coker(f)$ is a Gorenstein projective left $B$-module and $\coker(g)$ is a Gorenstein projective left $A$-module, then  $(X,Y,f,g)$ is a Gorenstein projective left $\Lambda$-module.\item If $(A,B,0,0)$ is a generalized compatible $\Lambda$-$\Lambda$-bimodule and $(X,Y,f,g)$ is a Gorenstein projective left $\Lambda$-module, then the sequences $V\otimes_B U\otimes_A X\stackrel{V\otimes f}\rightarrow V\otimes_B Y\stackrel{g}\rightarrow X$ and $U\otimes_A V\otimes_B Y\stackrel{U\otimes g}\rightarrow U\otimes_A X\stackrel{f}\rightarrow Y$ are exact, $\coker(f)$ is a Gorenstein projective left $B$-module and $\coker(g)$ is a Gorenstein projective left $A$-module.
 \end{enumerate}
\end{thm}
\begin{proof}(1) By Lemma \ref{lem: 5.1}(1), $U\oplus V$ is a generalized compatible $A\times B$-$A\times B$-bimodule. Since the sequences $V\otimes_B U\otimes_A X\stackrel{V\otimes f}\rightarrow V\otimes_B Y\stackrel{g}\rightarrow X$ and $U\otimes_A V\otimes_B Y\stackrel{U\otimes g}\rightarrow U\otimes_A X\stackrel{f}\rightarrow Y$ are exact, the sequence $(U\oplus V)\otimes_{A\times B} (U\oplus V)\otimes_{A\times B} (X,Y)\stackrel{(U\oplus V)\otimes(g,f)}\rightarrow (U\oplus V)\otimes_{A\times B} (X,Y)\stackrel{(g,f)}\rightarrow (X,Y)$ is exact. Since $\coker(f)$ is a Gorenstein projective left $B$-module and $\coker(g)$ is a Gorenstein projective left $A$-module, $\coker(g,f)=(\coker(g),\coker(f))$ is a Gorenstein projective left $A\times B$-module. By Theorem \ref{thm: 3.3}, $((X,Y),(g,f))$ is a Gorenstein projective left $(A\times B) \ltimes (U\oplus V)$-module. So $(X,Y,f,g)$ is a Gorenstein projective left $\Lambda$-module.

(2) Since $(A,B,0,0)$ is a generalized compatible $\Lambda$-$\Lambda$-bimodule, $\textbf{Z}(A\times B)$ is a generalized  compatible $(A\times B) \ltimes (U\oplus V)$-$(A\times B) \ltimes (U\oplus V)$-bimodule. Since $(X,Y,f,g)$ is a Gorenstein projective left $\Lambda$-module, $((X,Y),(g,f))$ is a Gorenstein projective left $(A\times B) \ltimes (U\oplus V)$-module. By Theorem \ref{thm: 3.4}, the sequence $(U\oplus V)\otimes_{A\times B} (U\oplus V)\otimes_{A\times B} (X,Y)\stackrel{(U\oplus V)\otimes(g,f)}\rightarrow (U\oplus V)\otimes_{A\times B} (X,Y)\stackrel{(g,f)}\rightarrow (X,Y)$ is exact and $\coker(g,f)$ is a Gorenstein projective left $A\times B$-module. So the sequences $V\otimes_B U\otimes_A X\stackrel{V\otimes f}\rightarrow V\otimes_B Y\stackrel{g}\rightarrow X$ and $U\otimes_A V\otimes_B Y\stackrel{U\otimes g}\rightarrow U\otimes_A X\stackrel{f}\rightarrow Y$ are exact, $\coker(f)$ is a Gorenstein projective left $B$-module and $\coker(g)$ is a Gorenstein projective left $A$-module.
\end{proof}
\begin{thm} \label{thm: 5.3}Let $\Lambda=\biggl(\begin{matrix} A&_{A}V_{B}\\_{B}U_{A}&B \end{matrix}\biggr)_{(0,0)}$ and $[X,Y,f,g]$ be a left $\Lambda$-module.
 \begin{enumerate}\item If both $_{B}U_{A}$ and $_{A}V_{B}$ are generalized cocompatible bimodules, the sequences $X\stackrel{f}\rightarrow \Hom_{B}(U,Y)\stackrel{{\rm Hom}_{B}(U,g)}\rightarrow \Hom_{B}(U,\Hom_{A}(V,X))$ and $Y\stackrel{g}\rightarrow \Hom_{A}(V,X)\stackrel{{\rm Hom}_{A}(V,f)}\rightarrow \Hom_{A}(V,\Hom_{B}(U,Y))$  are exact, $\ker(f)$ is a Gorenstein injective left $A$-module and $\ker(g)$ is a Gorenstein injective left $B$-module, then  $[X,Y,f,g]$ is a Gorenstein injective left $\Lambda$-module.\item If  $[A,B,0,0]$  is a generalized cocompatible $\Lambda$-$\Lambda$-bimodule and $[X,Y,f,g]$ is a Gorenstein injective left $\Lambda$-module, then the sequences $X\stackrel{f}\rightarrow \Hom_{B}(U,Y)\stackrel{{\rm Hom}_{B}(U,g)}\rightarrow \Hom_{B}(U,\Hom_{A}(V,X))$ and $Y\stackrel{g}\rightarrow \Hom_{A}(V,X)\stackrel{{\rm Hom}_{A}(V,f)}\rightarrow \Hom_{A}(V,\Hom_{B}(U,Y))$  are exact, $\ker(f)$ is a Gorenstein injective left $A$-module and $\ker(g)$ is a Gorenstein injective left $B$-module.
 \end{enumerate}
\end{thm}
\begin{proof}(1) By Lemma \ref{lem: 5.1}(2), $U\oplus V$ is a generalized cocompatible $A\times B$-$A\times B$-bimodule.

Since the sequences  $X\stackrel{f}\rightarrow \Hom_{B}(U,Y)\stackrel{{\rm Hom}_{B}(U,g)}\rightarrow \Hom_{B}(U,\Hom_{A}(V,X))$ and $Y\stackrel{g}\rightarrow \Hom_{A}(V,X)\stackrel{{\rm Hom}_{A}(V,f)}\rightarrow \Hom_{A}(V,\Hom_{B}(U,Y))$ are exact, it follows that the sequence $(X,Y)\stackrel{(f,g)}\rightarrow \Hom_{A\times B}(U\oplus V, (X,Y))\stackrel{(f,g)_{*}}\rightarrow\Hom_{A\times B}(U\oplus V,\Hom_{A\times B}(U\oplus V, (X,Y)))$ is also exact. Since $\ker(f)$ is a Gorenstein injective left $A$-module and $\ker(g)$ is a Gorenstein injective left $B$-module, $\ker(f,g)=(\ker(f),\ker(g))$ is a Gorenstein injective left $A\times B$-module. By Theorem \ref{thm: 4.3}, $[(X,Y),(f,g)]$ is a Gorenstein injective left $(A\times B) \ltimes (U\oplus V)$-module. So $[X,Y,f,g]$ is a Gorenstein injective left $\Lambda$-module.

(2) Since $[A,B,0,0]$  is a generalized cocompatible $\Lambda$-$\Lambda$-bimodule, $\textbf{Z}(A\times B)$ is a generalized  cocompatible $(A\times B) \ltimes (U\oplus V)$-$(A\times B) \ltimes (U\oplus V)$-bimodule. Since $[X,Y,f,g]$  is a Gorenstein injective left $\Lambda$-module, $[(X,Y),(f,g)]$  is a Gorenstein injective left $(A\times B) \ltimes (U\oplus V)$-module. By Theorem \ref{thm: 4.4}, the sequence  $(X,Y)\stackrel{(f,g)}\rightarrow \Hom_{A\times B}(U\oplus V, (X,Y))\stackrel{(f,g)_{*}}\rightarrow\Hom_{A\times B}(U\oplus V,\Hom_{A\times B}(U\oplus V, (X,Y)))$ is exact and $\ker(f,g)$ is a Gorenstein injective left $A\times B$-module. So the sequences  $X\stackrel{f}\rightarrow \Hom_{B}(U,Y)\stackrel{{\rm Hom}_{B}(U,g)}\rightarrow \Hom_{B}(U,\Hom_{A}(V,X))$ and $Y\stackrel{g}\rightarrow \Hom_{A}(V,X)\stackrel{{\rm Hom}_{A}(V,f)}\rightarrow \Hom_{A}(V,\Hom_{B}(U,Y))$  are exact, $\ker(f)$ is a Gorenstein injective left $A$-module and $\ker(g)$ is a Gorenstein injective left $B$-module.
\end{proof}
\begin{thm} \label{thm: 5.4}Let $\Lambda=\biggl(\begin{matrix} A&_{A}V_{B}\\_{B}U_{A}&B \end{matrix}\biggr)_{(0,0)}$  be a left coherent ring and $(W,Q,f,g)$ be a right $\Lambda$-module.
 \begin{enumerate}\item If  both $_{B}U_{A}$ and $_{A}V_{B}$ are generalized cocompatible bimodules, the sequences  $Q\otimes_B U\otimes_A V\stackrel{f\otimes V}\rightarrow W\otimes_A V\stackrel{g}\rightarrow Q$ and $W\otimes_A V\otimes_B U\stackrel{g\otimes U}\rightarrow Q\otimes_B U\stackrel{f}\rightarrow W$ are exact,  $\coker(f)$ is a Gorenstein flat right $A$-module and $\coker(g)$ is a Gorenstein flat right $B$-module,  then $(W,Q,f,g)$  is a Gorenstein flat right $\Lambda$-module.\item If  $[A,B,0,0]$  is a generalized   cocompatible $\Lambda$-$\Lambda$-bimodule  and  $(W,Q,f,g)$  is a  Gorenstein flat right $\Lambda$-module, then the sequences  $Q\otimes_B U\otimes_A V\stackrel{f\otimes V}\rightarrow W\otimes_A V\stackrel{g}\rightarrow Q$ and $W\otimes_A V\otimes_B U\stackrel{g\otimes U}\rightarrow Q\otimes_B U\stackrel{f}\rightarrow W$ are exact,  $\coker(f)$ is a Gorenstein flat right $A$-module and $\coker(g)$ is a Gorenstein flat right $B$-module.
 \end{enumerate}
\end{thm}
\begin{proof}(1) By Lemma \ref{lem: 5.1}(2), $U\oplus V$ is a generalized cocompatible $A\times B$-$A\times B$-bimodule. Since the sequences  $Q\otimes_B U\otimes_A V\stackrel{f\otimes V}\rightarrow W\otimes_A V\stackrel{g}\rightarrow Q$ and $W\otimes_A V\otimes_B U\stackrel{g\otimes U}\rightarrow Q\otimes_B U\stackrel{f}\rightarrow W$  are exact, the sequence $(W,Q)\otimes_{A\times B}(U\oplus V)\otimes_{A\times B} (U\oplus V)\stackrel{(f,g)\otimes(U\oplus V)}\rightarrow (W,Q)\otimes_{A\times B}(U\oplus V)\stackrel{(f,g)}\rightarrow (W,Q)$ is exact. Since  $\coker(f)$ is a Gorenstein flat right $A$-module and $\coker(g)$ is a Gorenstein flat right $B$-module,  $\coker(f,g)=(\coker(f),\coker(g))$ is a Gorenstein flat right $A\times B$-module. By Theorem \ref{thm: 4.7}(1), $((W,Q),(f,g))$ is a Gorenstein flat right $(A\times B) \ltimes (U\oplus V)$-module. So $(W,Q,f,g)$  is a Gorenstein flat right $\Lambda$-module.

(2) Since $[A,B,0,0]$ is a generalized cocompatible $\Lambda$-$\Lambda$-bimodule, $\textbf{Z}(A\times B)$ is a generalized  cocompatible $(A\times B) \ltimes (U\oplus V)$-$(A\times B) \ltimes (U\oplus V)$-bimodule. Since  $(W,Q,f,g)$  is a  Gorenstein flat right $\Lambda$-module,  $((W,Q),(f,g))$ is a Gorenstein flat right $(A\times B) \ltimes (U\oplus V)$-module. By Theorem \ref{thm: 4.7}(2), the sequence $(W,Q)\otimes_{A\times B}(U\oplus V)\otimes_{A\times B} (U\oplus V)\stackrel{(f,g)\otimes(U\oplus V)}\rightarrow (W,Q)\otimes_{A\times B}(U\oplus V)\stackrel{(f,g)}\rightarrow (W,Q)$ is exact and $\coker(f,g)$ is a Gorenstein flat right $A\times B$-module. So the sequences $Q\otimes_B U\otimes_A V\stackrel{f\otimes V}\rightarrow W\otimes_A V\stackrel{g}\rightarrow Q$ and $W\otimes_A V\otimes_B U\stackrel{g\otimes U}\rightarrow Q\otimes_B U\stackrel{f}\rightarrow W$ are exact,  $\coker(f)$ is a Gorenstein flat right $A$-module and $\coker(g)$ is a Gorenstein flat right $B$-module.
\end{proof}
\bigskip
\centerline {\bf ACKNOWLEDGEMENTS}
\bigskip
This research was supported by NSFC (12171230, 12271249) and NSF of Jiangsu Province of China (BK20211358).

\end{document}